\documentclass[12pt]{amsart}
\usepackage{graphicx, amsmath, amssymb, amsthm, amsfonts, mathtools, tikz - cd, url, comment, tikz, caption, subcaption} 

\newtheorem{thm}{Theorem}[section]
\newtheorem{lem}[thm]{Lemma}
\newtheorem{prop}[thm]{Proposition}
\newtheorem{cor}[thm]{Corollary}

\theoremstyle{definition}
\newtheorem{defn}[thm]{Definition}

\def\FS{\mathrm{FS}_\mathbb{N}}
\def\S{\mathrm{S}_\mathbb{N}}
\def\PS{\mathrm{PS}_\mathbb{N}}

\def\co{\mathrm{co}}
\def\oco{\overline{\mathrm{co}}}
\def\aff{\mathrm{aff}}
\def\oaff{\overline{\mathrm{aff}}}
\def\bmid{\text{ }\biggr\lvert\text{ }}
\def\DS{\mathrm{DS}}
\def\DSS{\mathrm{DSS}}
\def\ext{\mathrm{ext}}

\newenvironment{proofm}{%
\proof}{\endproof}

\title{An Operator Theoretic Approach to Birkhoff's Problem 111}
\author{Miles Gould}
\date{December 2024}

\begin{document}

\begin{abstract}
In 1946, Garrett Birkhoff proved that the $n\times n$ doubly stochastic matrices comprise the convex hull of the $n\times n$ permutation matrices, which in turn make up the extreme points of this polytope. He proposed his problem 111, which asks whether there exists a topology on infinite matrices for which this applies to the closed convex hull of the $\mathbb{N}\times\mathbb{N}$ permutation matrices. As Isbell showed in 1955, this equality is not achieved in the line-sum norm. In this paper, we use the domain of operator theory, and its many topologies, to improve on his negative result by showing that Birkhoff's problem is not solved in any of these topologies. In Kendall's 1960 paper on this problem, he gave an answer to the affirmative, as well as a topology for which closed convex hull comprises the doubly substochastic matrices. We also show that Kendall's secondary theorem also applies for all the locally convex Hausdorff topologies finer than than Kendall's (namely that of entry-wise convergence) which make the continuous dual of the matrix space no larger than the predual of the von Neumann algebra containing them. We then show that this is a theoretical upper limit topologies with this closure property. We also discuss the exposed points of this hull for these several topologies. Moreover, we show that, in these topologies, the closed affine hull of these permutation matrices comprise all operators with real-entry matrix coefficients.
\end{abstract}

\maketitle

\section{The Standard Representation of $\S$}

\begin{defn}
Let $H$ be a separable infinite dimensional Hilbert space. For any orthonormal basis $E=(e_n)_{n\in\mathbb{N}}$ of $H,$ we define the standard unitary representation $\pi^E:\S\rightarrow B(H)$ (or simply $\pi$ when $E$ is clear) of $\S$ by $\pi(\rho)e_n=e_{\rho(n)}.$
\end{defn}

For the reader strictly interested in our considerations on convex hulls, this entire first section need not be read, as it serves only to justify the use of the standard representation $\pi.$

Our definition of a unitary representation $\pi$ of a topological group $G$ is that of \cite{BE19}, which in addition to being a group homomorphism into the unitaries $U(H),$ must make also the map $G\times H\rightarrow H$ by $(g,x)\mapsto \pi(g)x$ continuous for the product topology. For any of the mainstream operator topologies, namely weak, ultraweak, strong, strong$^*$, ultrastrong, ultrastrong$^*$, Arens-Mackey, weak Banach, and norm, this map is continuous. Per \cite{BE19}, the continuity of this map is equivalent to the continuity of $\pi:G\rightarrow U(H)$ where the codomain is equipped with the strong topology. Moreover, we know the weak and strong coincide on $U(H),$ so for any topology finer than the weak, this continuity follows.

\begin{prop}\label{pwtop}
The pointwise topology is the coarsest topology which makes $S_\mathbb{N}$ a Hausdorff topological group. Moreover, it is the coarsest topology on $S_\mathbb{N}$ which makes $\pi$ continuous with respect to the weak topology on $B(H).$
\end{prop}

\begin{proofm}
The first part of the theorem is given by \cite{GA67}. Next, a net $(\pi_\lambda)$ in $B(H)$ converges to $\pi$ strongly in $B(H)$ if and only if, for all $x\in H,$
\[\lim_\lambda\sum_{n=1}^\infty|\langle (\pi_\lambda-\pi)x,e_n\rangle|^2=0.\]
Therefore, if $S_\mathbb{N}$ is equipped with the topology we seek, a net $(\rho_\lambda)$ converges to $\rho$ if and only if
\[\lim_\lambda\sum_{n=1}^\infty|\langle (\pi_\lambda(\rho)-\pi(\rho))x,e_n\rangle|^2=0,\]
which simplifies to
\[\lim_\lambda\sum_{n=1}^\infty|x_{\rho^{-1}_\lambda(n)}-x_{\rho^{-1}(n)}|^2=0.\]
Notice the following bound
\begin{align*}
\lim_\lambda\sum_{n=1}^\infty|x_{\rho^{-1}_\lambda(n)}-x_{\rho^{-1}(n)}|^2 &\leq \sum_{n=1}^\infty\left(|x_{\rho^{-1}_\lambda(n)}|^2+|x_{\rho^{-1}(n)}|^2\right) \\
&\leq \sum_{n=1}^\infty|x_{\rho^{-1}_\lambda(n)}|^2+\sum_{n=1}^\infty|x_{\rho^{-1}(n)}|^2 \\
&\leq 2\sum_{n=1}^\infty|x_n|^2 \\
&\leq 2\|x\| \\
&<\infty.
\end{align*}
Thus, the prior limit holds if and only if, for each $n\in\mathbb{N},$ $x_{\rho^{-1}_\lambda(n)}\rightarrow x_{\rho^{-1}(n)}.$ Moreover, since this must hold for every $x\in H,$ convergence holds if and only if $\rho^{-1}_\lambda(n)\rightarrow\rho^{-1}(n).$ Since $\mathbb{N}$ is discrete, this holds if and only if $\rho^{-1}_\lambda(n)=\rho^{-1}(n)$ for eventual $\lambda.$ Finally, because every $n=\rho(m)$ for some unique $m\in\mathbb{N},$ convergence is equivalent to the statement: for all $m\in\mathbb{N},$ $\rho_\lambda(m)=\rho(m)$ for eventual $\lambda.$
\end{proofm}

An interesting consequence of equipping $S_\mathbb{N}$ with the pointwise topology is that every irreducible unitary representation of $S_\mathbb{N}$ must be over a separable Hilbert space. Therefore the degree of every faithful irreducible unitary representation is bounded above by $\aleph_0.$ We will now show that the degree must also be bounded below by $\aleph_0.$

\begin{lem}
Let $n\in\mathbb{N}.$ Then every unitary representation $\pi:S_\mathbb{N}\rightarrow M_n(\mathbb{C})$ is unfaithful.
\end{lem}

\begin{proofm}
Suppose, towards contradiction that $\pi:S_\mathbb{N}\rightarrow M_n(\mathbb{C})$ is faithful. Let $\sim$ be the equivalence relation on $S_\mathbb{N}$ by $\rho\sim\nu$ iff $\rho\nu=\nu\rho.$ By assumption, the quotient $S_\mathbb{N}/\sim$ must be finite. Consider the sequence $(\rho_k)$ in $S_\mathbb{N}$ by $\rho_k:(1\text{ }2\cdots k+1)$ in cycle notation. Whenever $m\neq k,$ $\rho_k\rho_m\neq \rho_m\rho_k,$ a contradiction, since this sequence could not fit into $S_\mathbb{N}/\sim.$
\end{proofm}

\begin{prop}
For any topology $\mathcal{T}$ on $S_\mathbb{N}$ finer than pointwise, $\pi$ is a faithful non-degenerate unitary representation of $(S_\mathbb{N},\mathcal{T})$ on $H.$
\end{prop}

\begin{proofm}
Let $\rho,\nu\in S_\mathbb{N}.$ Then
\[\pi(\rho)\pi(\nu)(e_n)=\pi(\rho)(e_{\nu(n)})=e_{\rho\nu(n)},\]
so $\pi$ is a group homomorphism (and so a representation). Moreover, 
\[\langle\pi(\rho)(e_n),e_m\rangle=\delta_{\rho(n),m}\]
where $\delta$ is the Kronecker delta. Since it is real, 
\[\delta_{\rho(n),m}=\delta_{\rho^{-1}(m),n}=\langle \pi(\rho)^*(e_m),e_n\rangle,\]
so $\pi(\rho)^*=\pi(\rho^{-1})=\pi(\rho)^{-1}.$ Therefore $\pi$ is unitary. Now suppose $\pi(\rho)=\pi(\nu).$ Then $e_{\rho(n)}=e_{\nu(n)}$ for all $n\in\mathbb{N},$ i.e. $\rho=\nu,$ so $\pi$ is faithful. Since $\mathrm{Id}_H\in\pi(S_\mathbb{N}),$ non-degeneracy is trivial.
\end{proofm}

Next, we will show a property which differs between the standard representation of the infinite symmetric group and the standard representations of the finite symmetric groups, namely that the former is irreducible.

\begin{lem}\label{finrep}
Let $\pi_m: S_m\rightarrow M_m(\mathbb{C})$ by $\pi_m(\rho)e_n=e_{\rho(n)}.$ Then the commutant $\pi(S_m)'=\mathrm{span}(\mathrm{Id}_m,U_m),$ where $[U_m]_{ij}=1.$ In particular, $\pi_m$ is reducible.
\end{lem}

\begin{proofm}
(Induction) Clearly, $U_1=\mathrm{Id}_1,$ so $\pi_1(S_1)'=\mathbb{C}\mathrm{Id}_1$ is trivial. Now, let $k\in\mathbb{N}$ and suppose $\pi(S_k)'=\mathrm{span}(\mathrm{Id}_k,U_k).$ Let 
$a\in\pi_{k+1}(S_{k+1}).$ We can express it in $k\times k$ and $1\times 1$ block form. Then for every $\rho\in S_k,$

\begin{align*}
\begin{pmatrix}
\pi_k(\rho) & 0 \\
0 & 1
\end{pmatrix}
\begin{pmatrix}
A & B \\
C & D
\end{pmatrix}
&=
\begin{pmatrix}
A & B \\
C & D
\end{pmatrix}
\begin{pmatrix}
\pi_k(\rho) & 0 \\
0 & 1
\end{pmatrix} \\
\begin{pmatrix}
\pi_k(\rho)A & \pi_k(\rho)B \\
C & D
\end{pmatrix}
&=
\begin{pmatrix}
A\pi_k(\rho) & B \\
C\pi_k(\rho) & D
\end{pmatrix}
\end{align*}
Then $A\in \pi_k(S_k)',$ $B$ is a $1\times k$ constant row, $C$ a $k\times 1$ constant column, and $D$ is free.
We can then do the same for the $1\times 1$ and $k\times k$ block form, yielding the general equation
\begin{align*}
\begin{pmatrix}
a_1 & a_2 & \cdots & a_2 & b \\
a_2 & a_1 & \cdots & a_2 & b \\
\vdots & \vdots & \ddots & \vdots & \vdots \\
a_2 & a_2 & \cdots & a_1 & b \\
c & c & \cdots & c & d
\end{pmatrix}
&=
\begin{pmatrix}
a' & b' & \cdots & b' & b' \\
c' & d'_1 & \cdots & d'_2 & d'_2 \\
\vdots & \vdots & \ddots & \vdots & \vdots \\
c' & d'_2 & \cdots & d'_1 & d'_2 \\
c' & d'_2 & \cdots & d'_2 & d'_1
\end{pmatrix} \\
a'=a_1 &= d=d'_1 \\
a_2=d'_2 = b' &= b=c'=c
\end{align*}
and so we have the final form $a\in\mathrm{span}(\mathrm{Id}_m,U_m).$ It is simple to check that all such matrices commute with $\pi_{k+1}(S_{k+1})$. By induction, the lemma holds.
\end{proofm}

\begin{thm}
For any topology finer than pointwise, $\pi$ is irreducible.
\end{thm}

\begin{proofm}
As before, we can express a general element $a\in\pi(S_\mathbb{N})'$ in $m\times m$ and $\mathbb{N}\times\mathbb{N}$ block form, and we can see that, for all such block decompositions $A\in\mathrm{span}(\mathrm{Id}_m,U_m),$ so every matrix coefficient must fit this pattern. However, if $a\not\in \mathbb{C}\mathrm{Id}_H$ then $a$ must be unbounded, a contradiction. Thus, $\pi(S_\mathbb{N})'=\mathbb{C}\mathrm{id}_H.$
\end{proofm}

\section{The Convex Hulls of $\S$}

Quickly, we will show that the group algebra of $\S$ on $H$ is strongly dense in $B(H),$ which we will analogize for real-entry matrices later.

\begin{cor}\label{dense}
$\mathbb{C}[\pi(S_\mathbb{N})]$ is strongly dense in $B(H).$ Moreover, 
\\$\mathbb{C}[\pi(S_\mathbb{N})]\cap U(H)$ is strongly dense in $U(H).$
\end{cor}

\begin{proofm}
This follows from the double commutant lemma and Kaplansky's density theorem \cite[123,131]{MU90}.
\end{proofm}

\begin{defn}
For subsets $S\subseteq B(H),$ we will use the following notation for particular closures of $S,$ $\overline{S}$ the uniform closure, $\overline{S}^s$ the strong closure, $\overline{S}^w$ the weak, etc.
\end{defn}

We are now going to define several crucial collections of permutations, partial permutations, and their polytopes.

\begin{defn}
Firstly, we should note that, just as $\pi$ depends on $E,$ so will all the following definitions, via the matrix coefficients $u_{mk}=\langle ue_k,e_m\rangle.$
\\ (i) $\DS=\biggr\{a\in B(H)\bmid a_{mk}\geq 0\text{, }\sum_{k=1}^\infty a_{mk}=\sum_{k=1}^\infty a_{km}=1\biggr\},$
the doubly stochastic matrices.
\\ (ii) $\DSS=\biggr\{a\in B(H)\bmid a_{mk}\geq 0\text{, }\sum_{k=1}^\infty a_{mk}\leq 1\text{, }\sum_{k=1}^\infty a_{km}\leq 1\biggr\},$
the doubly substochastic matrices.
\\ (iii) $\PS=\biggr\{a\in B(H)\bmid a_{mk}\in\{0,1\}\text{, }\sum_k a_{mk}\leq 1\text{, }\sum_k a_{km}\leq 1\biggr\},$ the partial permutation matrices. Their finite analog has been studied recently, by \cite{ST12}.
\\ (iv) $\FS=\{\rho\in \S\mid \rho(n)=n\text{ e.a.}\}$ the finitary symmetric group.
\end{defn}

\subsection{The Uniform and Finitary Convex Hulls}

\begin{prop}\label{normco}
$\oco^{wB}\pi(\S)\subset \DS$
\end{prop}

\begin{proofm}
Following in the footsteps of Isbell \cite{IS55}, we will show by example that the above inclusion is strict. Let $a\in\DS$ by
\[a=\begin{pmatrix}
1 & 0 & 0 & 0 & 0 & 0 & \cdots \\
0 & 1/2 & 1/2 & 0 & 0 & 0 & \cdots \\
0 & 1/2 & 1/2 & 0 & 0 & 0 & \cdots \\
0 & 0 & 0 & 1/3 & 1/3 & 1/3 & \cdots \\
0 & 0 & 0 & 1/3 & 1/3 & 1/3 & \cdots \\
0 & 0 & 0 & 1/3 & 1/3 & 1/3 & \cdots \\
\vdots & \vdots & \vdots & \vdots & \vdots & \vdots &\ddots \\
\end{pmatrix}\]
Let $b\in\co(\S),$ and express it $b=\sum_{j=1}^pt_j\pi(\rho_j).$ Then consider the $n\times n$ block of $1/n$ entries. Then, there are at most $p^2$ columns in this block for which $b$ has a nonzero entry. Take $x_k=\frac{1}{\sqrt{n-p^2}}$ for each of the remaining columns. 
We can lower bound the norm difference
\begin{align*}
\|a-b\|^2 &= \sup_{\|x\|=1}\sum_{m=1}^\infty\biggr\lvert\sum_{k=1}^\infty(a_{mk}-b_{mk})x_k\biggr\lvert^2 \\
&\geq n\biggr\lvert\frac{n-p^2}{n}\frac{1}{\sqrt{n-p^2}}\biggr\lvert^2 \\
&\geq \frac{n-p^2}{n}
\end{align*}
Since we can take $n$ arbitrarily large, we can ensure $\|a-b\|$ arbitrarily close to 1, so $\|a-b\|\geq 1$. In fact,
\[\inf_{b\in \co(\S)}\|a-b\|\geq1.\]
Now, since the weak Banach dual coincides with the norm dual, by \cite[370]{TR67}, $\oco^{wB}\pi(\S)\subset \DS.$
\end{proofm}

Our further work will show that $\oco^\mathcal{T}\pi(\S)\neq \DS$ for any of the 9 main operator topologies, thus yielding a negative answer to Birkhoff's problem 111 for operator topologies. Crucial to the strong and weak topologies, we will define two types of finite restrictions of operators with respect to $E$.

\begin{defn}
For $u\in B(H),$ let $u^{[n]},u^{\langle n\rangle}\in B(H)$ by
\[u^{\langle n\rangle}_{mk} = \begin{cases} 
        u_{mk}, & m,k\leq n; \\
        0, & \text{otherwise}.
   \end{cases}\]
and
\[u^{[n]}_{mk} = \begin{cases} 
        u_{mk}, & m\leq n\text{ or }k\leq n; \\
        0, & \text{otherwise}.
   \end{cases}\]

Now, for a subset $S\subseteq B(H),$ let
\[F_0S=\{v\in B(H)\mid \exists u\in S,n\in\mathbb{N}\text{ s.t. }u^{\langle n\rangle}=v\}\]
\[F_1S=\{v\in B(H)\mid \exists u\in S,n\in\mathbb{N}\text{ s.t. }u^{[n]}=v\}\]
We will also define the finitary $\textit{restriction}$ $FS$ by
\[FS=\{v\in S\mid \exists u\in S\text{, }n\in\mathbb{N}\text{ s.t. }a=b^{\langle n\rangle}+\mathrm{id}_{n+1}\},\]
where $\mathrm{id}_n\in B(H)$ by
\[\mathrm{id}_ne_k = \begin{cases} 
        0, & k\leq n; \\
        e_k, & k>n.
   \end{cases}\]
Notice that it is possible for $FS=\emptyset,$ namely if every $a\in S$ is never eventually the identity. However, some of our crucial examples behave nicely under $F,$ as shown the next lemma.
\end{defn}

\begin{lem}\label{finex}
$F\pi(\S)=\pi(\FS)$ and $F\DS=\bigcup_{n\in\mathbb{N}}(B_n+\mathrm{id}_{n+1}),$
where $B_n$ is the Birkhoff polytope identified in the upper-left $n\times n$ subalgebra.
\end{lem}

\begin{proofm}
If for $u\in\pi(\S),$ $u^{\langle n\rangle}+\mathrm{id}_{n+1}\in \pi(\FS),$ then necessarily $u^{\langle n\rangle}$ permutes $\{e_1,...,e_n\},$ so $u^{\langle n\rangle}+\mathrm{id}_{n+1}\in\pi(\FS).$ In the other direction, if $u\in\pi(\FS),$ then for some $n\in\mathbb{N},$ $u=u^{\langle n\rangle}+\mathrm{id}_{n+1}.$

Now if $u\in\DS$ has $u^{\langle n\rangle}+\mathrm{id}_{n+1}\in\DS,$ then necessarily $u^{\langle n\rangle}$ is doubly stochastic in the $n\times n$ subalgebra, so is in $B_n.$ The other direction is trivial.
\end{proofm}

\begin{prop}
$\co\pi(\FS)=F\DS.$
\end{prop}

\begin{proofm}
By lemma $\ref{finex},$ $\co\pi(\FS)=\co(F\pi(\S)).$ Now, due the fact that $B_n\subseteq B_{n+1}$ and a quick application of the Birkhoff-von Neumann theorem, $\co(F\pi(\S))=\bigcup_{n\in\mathbb{N}}(B_n+\mathrm{id}_{n+1})=F\DS.$
\end{proofm}

\begin{lem}\label{fincl}
For all $S\subseteq B(H),$ $S\subseteq \overline{F_0S}^w$ and $S\subseteq \overline{F_1S}^s.$
\end{lem}

\begin{proofm}
Let $u\in S,$ $x\in H.$ By Parseval's identity,
\begin{align*}
\|(u-u^{[n]})x\|^2 &= \sum_{k=1}^\infty\biggr\lvert\sum_{m=1}^\infty(u_{km}-u^{[n]}_{km})x_m\biggr\lvert^2 \\
&= \sum_{k=n+1}^\infty\biggr\lvert\sum_{m=n+1}^\infty u_{km}x_m\biggr\lvert^2 \\
&\rightarrow 0.
\end{align*}
Additionally, let $y\in H.$
\begin{align*}
|\langle (u-u^{\langle n\rangle})x,y\rangle| &= \biggr\lvert\sum_{m=1}^\infty\sum_{k=1}^\infty(u_{mk}-u^{\langle n\rangle}_{mk})x_k\overline{y_m}\biggr\lvert \\
&= \biggr\lvert\sum_{m=1}^n\sum_{k=n+1}^\infty u_{mk}x_k\overline{y_m}+\sum_{m=n+1}^\infty\sum_{k=1}^\infty u_{mk}x_k\overline{y_m}\biggr\lvert \\
&= \biggr\lvert\sum_{m=1}^n\sum_{k=1}^\infty u_{mk+n}x^{(n)}_k\overline{y_m} +\sum_{m=1}^\infty\sum_{k=1}^\infty u_{m+nk}x_k\overline{y^{(n)}_m}\biggr\lvert \\
&\leq \|u\|\left(\|x^{(n)}\|\|y\|+\|x\|\|y^{(n)}\|\right) \\
&\rightarrow 0
\end{align*}
Here we used the notation $x^{(n)}$ for the vector by 
\[x^{(n)}_k = \begin{cases} 
        0, & k\leq n; \\
        x_{n+k}, & n > k.
   \end{cases}\]
\end{proofm}

\subsection{The Arens-Mackey Family of Hulls}

Here we will consider a crucial family of operator topologies, namely those which are coarser than the Mackey topology yet finer than the weak topology. We will begin with the closures of the finitary permutation matrices in the strong and weak topologies.

\begin{prop}\label{clos}
$\overline{\pi(\FS)}^s=\pi(\S)$ and
$\overline{\pi(\FS)}^w=\PS.$
\end{prop}

\begin{proofm}
By cor. $\ref{fincl},$ $\pi(\S)\subseteq \overline{\pi(\FS)}^s,$ because every permutation matrix $u\in\pi(\S)$ has convergent $u^{[n]}$ (in fact it will coincide with $u^{\langle n\rangle}).$

Let $u\in\overline{\pi(\FS)}^w.$ Then there exists a sequence in $\pi(S_\mathbb{N})$ all $m,k\in\mathbb{N},$ $\lim_{n\rightarrow\infty}\langle(\pi(\rho_n)-u)e_k,e_m\rangle=0,$ but $\langle\pi(\rho_n)e_k,e_m\rangle=\delta_{m\rho_n(k)}.$ Clearly then, $u_{mk}\in \{0,1\}.$ Suppose for some $m,k,j,$
$u_{mk}=u_{mj}=1.$ Then
\begin{align*}
\lim_{n\rightarrow\infty}\delta_{m\rho_n(k)} &= \lim_{n\rightarrow\infty}\delta_{m\rho_n(j)} \\
\rho_n(k) &= \rho_n(j) \text{ for e.a. } n\text{, so} \\
k &= j.
\end{align*}
The same argument can be applied to a fixed column, so each row and column must have at most one $1,$ i.e. $u\in\PS.$

By $\ref{fincl},$ we need only show $F_0\PS\subseteq \overline{\pi(\FS)}^w,$ so let $u\in F_0\PS.$ Let 
\[I=\{k\in\mathbb{N}\mid\forall m\text{, }u_{mk}=0\}\text{ and }J=\{m\in\mathbb{N}\mid \forall k\text{, }u_{mk}=0\}.\] 
Let $i:\mathbb{N}\rightarrow I,$ $j:\mathbb{N}\rightarrow J$ be their increasing enumerations. 
Define the finitary permutations $\rho_n$ by
\[\rho_n(k) = \begin{cases} 
        n+k, & k\leq n; \\
        n-k, & n< k\leq 2n; \\
        n, &\text{ otherwise.}

   \end{cases}\]

Define $u^{(n)}$ by
\[u^{(n)}_{mk} = \begin{cases} 
        u_{mk}, & m\in I\text{ or }k\in J; \\
        \pi(\rho_n)_{i(m),j(k)}, & \text{ otherwise}.

   \end{cases}\]
Finally, we need only show $\pi(\rho^{(n)})\rightarrow 0$ weakly. Similar to the proof of $\ref{fincl},$ let $x,y\in H.$
\begin{align*}
|\langle\pi(\rho^{(n)})x,y\rangle| &= \sum_{m=1}^\infty x_{\rho_n^{-1}(m)}\overline{y_m} \\
&= \sum_{m=1}^n(x_{n+m}\overline{y_m}+x_m\overline{y_{n+m}})+\sum_{m=2n+1}^\infty x_m\overline{y_m} \\
&\leq \|x^{(n)}\|\|y\|+\|x\|\|y^{(n)}\|+\|x^{2n}\|\|y^{2n}\| \\
&\rightarrow 0,
\end{align*}
maintaining the notation from $\ref{fincl}.$

For the strong case, since $\|(\pi(\rho_n)-\pi(\rho_m))x\|\rightarrow 0$, we can see that for fixed $k,$ $\rho_n(k)$ must converge, so every row and by the same argument, every column of $u$, must contain a $1.$ Therefore, they must contain exactly one $1,$ which is precisely the definition of $\pi(\S).$
\end{proofm}

The next lemma makes obvious the property which ensures the weak closure of $\DSS.$

\begin{lem}\label{finsum}
Let $u\in B(H).$ If for some $n,m\in\mathbb{N},$
\[\sum_{k=1}^n u_{mk}>1\text{ or }\sum_{k=1}^n u_{km}>1,\]
then $u\not\in \overline{\DSS}^w.$
\end{lem}

\begin{proofm}
Let $u,n,m$ be as above. Without loss of generality, we will assume that the $m^{th}$ row sums to greater than 1. Suppose, towards contradiction, that $(u^{(j)})$ in $\PS$ converges to $u.$ Then for all $k\leq n,$
\[\langle (u^{(j)}-u)e_k,e_m\rangle\rightarrow 0.\]
Therefore for eventual $j,$
\[\sum_{k=1}^nu^{(j)}_{mk}>1,\]
a contradiction to $\DSS.$
\end{proofm}

This lemma shows us the intuitive fact that doubly substochastic matrices lie within the unit ball of $B(H).$ Interestingly, this is true without any knowledge of the boundedness of the operators, only their matrix coefficients.

\begin{lem}\label{unit}
$\DSS\subseteq B(H)_1$
\end{lem}

\begin{proofm}
Using an argument shown to me by David Gao \cite{GA24}, we will first show that if $(k_i),(t_i)$ are nonnegative sequences such that $\sum_{i}k_i\leq 1,$ then $(\sum_ik_it_i)^2\leq \sum_{i}k_it_i^2.$ If $\sum_ik_i=0,$ then $k_i=0$ always, so this is obvious, otherwise,
\begin{align*}
\left(\sum_{i=1}^\infty k_it_i\right)^2 &= \left(\sum_{i=1}^\infty k_i\right)^2\left(\frac{\sum_{i=1}^\infty k_it_i}{\sum_{i=1}^\infty k_i}\right)^2 \\
&\leq \left(\sum_{i=1}^\infty k_i\right)^2\frac{\sum_{i=1}^\infty k_it_i^2}{\sum_{i=1}^\infty k_i} \text{ by Jensen's inequality}\\
&\leq \left(\sum_{i=1}^\infty k_it_i\right)\sum_{i=1}^\infty k_it_i^2 \\
&\leq \sum_{i=1}^\infty k_it_i^2.
\end{align*}
Let $x\in H$ s.t. $\|x\|=1,$ $a\in \DSS.$
\begin{align*}
\|ax\|^2 &\leq \sum_{m=1}^\infty\biggr\lvert\sum_{k=1}^\infty a_{mk}x_k\biggr\lvert^2 \\
&\leq \sum_{m=1}^\infty\left(\sum_{k=1}^\infty a_{mk}|x_k|\right)^2 \\
&\leq \sum_{m=1}^\infty \sum_{k=1}^\infty a_{mk}|x_k|^2 \\
&\leq \sum_{m=1}^\infty \left(\sum_{k=1}^\infty a_{mk}\right)|x_k|^2 \\
&\leq \sum_{m=1}^\infty|x_k|^2 = \|x\|^2=1.
\end{align*}
\end{proofm}

\begin{lem}\label{hullc}
$\DSS$ is convex and weakly compact. 
\end{lem}

\begin{proofm}
Let $u,w\in\DSS,$ $t\in(0,1).$ Then
\begin{align*}
\sum_{k=1}^\infty (tu_{mk}+(1-t)w_{mk})\leq 1
\end{align*}
so convexity holds. Since, by $\ref{unit},$ $\DSS\subseteq B(H)_1,$ which by Banach-Alaoglu is weakly compact, it suffices to show closure. Let $u\in B(H)$ such that there exists $m\in\mathbb{N}$ whose row or column has sum greater than 1 (i.e. $u\not\in\DSS)$. Without loss of generality, assume the row case. Then there exists some $n$ such that
\[\sum_{k=1}^nu_{mk}>1,\]
so lemma $\ref{finsum},$ $u$ is not in the weak closure of $\DSS.$ Therefore, $\DSS$ is closed.
\end{proofm}

\begin{lem}\label{extinc}
$\oco^w(\PS)\subseteq \DSS.$
\end{lem}

\begin{proofm}
By lemma $\ref{hullc},$ $\overline{co}^w(\DSS)=\DSS,$ so
\[\oco^w(\PS)\subseteq \oco^w(\DSS)=\DSS.\]
\end{proofm}

The following theorem is a strengthening of Kendall's result \cite{KE60}, who proved the result for the weakest topology which makes the functionals $a\mapsto a_{mk}$ continuous. We have already given (many) examples of sequences in $\DSS$ which converge in the weak and diverge in the Mackey topology, via $\ref{clos},$ and it is not difficult to show that the weak topology coincides with Kendall's in the unit ball.

\begin{thm}\label{coh}
$\overline{co}^M\pi(\FS)=\DSS$ and $\mathrm{ext}(\DSS)=\PS.$
\end{thm}

\begin{proofm}
Let $u\in\DSS,$ $n\in\mathbb{N}.$ The restriction $u^{\langle n\rangle}$ has
\[u^{\langle n\rangle}_{mk}\in[0,1]\text{, } \sum_{k=1}^nu^{\langle n\rangle}_{mk}\leq 1\text{, and }\sum_{k=1}^nu^{\langle n\rangle}_{km}\leq 1.\] 
Thus $u^{(n)}$ is in the partial permutohedron \cite{ST12}, by Mirsky \cite{MI59}. Therefore, \[u^{\langle n\rangle}\in \co(\PS^{\langle n\rangle})\subseteq\co(F_0\PS).\] 
Now, by lemma $\ref{fincl}$ and $\ref{clos},$ 
\[u\in \oco^w(F_0\PS)=\oco^w(\pi(\FS)).\]
The other inclusion is given by lemma $\ref{extinc}.$ Since $\DSS$ is contained in the unit ball, by $\ref{unit}$, the weak and ultraweak topologies coincide on it. Additionally, since by $\ref{hullc},$ $\DSS$ is convex, the ultraweak and Mackey closures are characterized by their functionals \cite[370]{TR67}, which coincide. Krein-Milman \cite[273]{MU90} then tells us that $\PS$ contains the extreme points of $\DSS.$ Finally, we can verify that every $v\in \PS$ is extreme. If $v=tu+(1-t)w$ for $u,w\in \DSS$ and $t\in(0,1),$ then for some $m,k,$ $v_{mk}\in[t,1-t]\subseteq(0,1),$ a contradiction.
\end{proofm}

It should be noted that we can achieve the above result in a briefer argument: As Kendall noted, $\ext(\DSS)=\PS,$ (he claimed this was easily calculated, though did not provide proof) so by lemma $\ref{hullc},$ Krein-Milman gives us the fact immediately. However, Krein-Milman requires the axiom of choice, so the proof outlined above is superior.

\begin{cor}
$\oco^\mathcal{T}\pi(S_\mathbb{N})\neq \DS$ for $\mathcal{T}$ in the following: weak, ultrweak, strong, strong$^*$, ultrastrong, ultrastrong$^*$, Arens-Mackey, weak Banach, norm. In particular, Birkhoff's problem 111 cannot be solved using operator topologies.
\end{cor}

\begin{proofm}
By prop. $\ref{normco},$ 
\[\oco\pi(\S)=\oco^{wB}\pi(\S)\subset \DS,\]
and by theorem $\ref{coh},$
\[\oco^w\pi(\S)=\oco^M\pi(\S)=\DSS\neq \DS,\]
so since the other topologies are finer than weak, though coarser than Mackey,
\[\oco^{uw}\pi(\S)=\oco^s\pi(\S)=\oco^{s^*}\pi(\S)=\oco^{us}\pi(\S)=\oco^{us^*}\pi(\S)=\DSS.\]
\end{proofm}

\subsection{Distinct Topologies which Generate $\mathbf{DSS}$}

The fact that the closed convex hulls of $\pi(\FS)$ coincide between the Mackey and weak topologies is actually quite remarkable (though not specific to $\pi(\FS)).$ As Kendall noted, $\DSS$ is compact in the entry-wise topology, so in the weak topology. However, theorem $\ref{coh}$ gives us two distinct topologies for which $\DSS$ is a non-compact closure, namely strong and strong$^*$. 

Let us summarize the existing arguments: It is well known that the weak and ultraweak topologies coincide on the unit ball, so by lemma $\ref{unit},$ the same is true of $\DSS.$ Moreover, as proven by von Neumann \cite{VO36}, we may say the same of the strong and ultrastrong. A bit later, Akemann \cite{AK67} proved that the strong$^*$ and Mackey topologies coincide on the unit ball. The next two propositions will show that these three classes are actually distinct on $\DSS.$

\begin{prop}
$\DSS$ is not strongly compact.
\end{prop}

\begin{proofm}
The proof follows the analogous result for the unit ball. Suppose, towards contradiction, that $\DSS$ is compact in the strong topology. Then the identity map $(\DSS,s)\mapsto (\DSS,w)$ is a continuous bijection from a compact space to a Hausdorff space, so is a homeomorphism. The prior continuity is ensured because the domain topology is finer than the codomain topology. However, in light of $\ref{clos},$ there are many weakly convergent sequences which diverge strongly, a contradiction.
\end{proofm}

This argument clearly holds for all the operator topologies other than the weak and ultraweak. This begs the question of how much finer (than Mackey) a topology can be so that theorem $\ref{coh}$ still holds. Interestingly, one can give a proof of the above using Krein-Milman: Suppose $\DSS$ is strongly compact. Then since $\pi(\S)$ is closed set and $\oco^s\pi(\S)=\DSS,$ then $\ext(\DSS)\subseteq\pi(\S),$ a contradiction. Again though, Krein-Milman relies on the axiom of choice, so the original proof is superior.

\begin{prop}
$(\DSS,s)$ and $(\DSS,s^*)$ are not homeomorphic.
\end{prop}

\begin{proofm}
We will prove this by example. Let $u^{(n)}\in\DSS$ by
\[u^{(n)}_{mk} = \begin{cases} 
        1, & m=1,k=n; \\
        0, & \text{otherwise}.

   \end{cases}\]
Let $x\in H.$ Then $\|u^{(n)}x\|=|x_n|\rightarrow0.$ However, $\|u^{(n)*}e_1\|=\|e_n\|=1,$ so $u^{(n)}$ converges strongly but not strongly$^*.$
\end{proofm}

Our next goal will be to prove that the strong$^*$ (or equivalently, Mackey) topology is the finest topology on $\DSS$ having $\co\pi(\FS)$ as a dense subset which is inherited from a locally convex topology on $B(H).$

\begin{defn}
Our further considerations on $\DSS$ require us to delve into the singular functionals, that is, functionals which do not correspond to the pre-dual. The mere existence of one of these demands the axiom of choice. Let $B(H)_s^*$ denote these singular functionals. That is,
\[B(H)_s^*=\mathrm{span}\{f\in B(H)_+^*\mid \forall g\in B(H)_{*+}\text{, }g\not\leq f\}.\]
\end{defn}

\begin{lem}\label{funcon}
If $f\in B(H)^*_s,$ then $f$ is constant on $\pi(\FS).$
\end{lem}

\begin{proofm}
Let $f\in B(H)^*_s.$ Let $u,w\in\co\pi(\FS).$ Then there exists $n\in\mathbb{N}$ s.t. $u-w=(u-w)^{\langle n\rangle}.$ Let $K_n=\mathrm{span}\{e_j\mid j\leq n\}.$ By Akemann \cite{AK67}, for all $x\in K_n,$ $f(x\otimes x)=0.$ Then $(u-w)^{\langle n\rangle}\in B(K_n),$ which is spanned by its rank-one projections, thus $f(u-w)=f((u-w)^{\langle n\rangle})=0,$ so $f(u)=f(w).$
\end{proofm}

\begin{lem}\label{funzero}
Let $\mathcal{T}$ be a locally convex topology on $B(H).$ If there exists a $\tau-$continuous functional $f\in B(H)_s^*$ such that $f(1)\neq 0,$ then $0\not\in\oco^\mathcal{T}\pi(\FS).$ 
\end{lem}

\begin{proofm}
By $\ref{funcon},$ if $f(1)\neq 0,$ then for all $u\in \co\pi(\FS),$ $f(u)=f(1)\neq 0,$ so no net in $\co\pi(\FS)$ approaches 0, otherwise $f$ could not be continuous, since $f(0)=0.$
\end{proofm}

\begin{cor}\label{closzero}
Let $\mathcal{T}$ be a locally convex topology on $B(H).$ If 
\[\oco^\mathcal{T}\pi(\FS)=\DSS,\]
then every $f\in B(H)^*_s$ is identically 0 on $\DSS.$
\end{cor}

\begin{proofm}
Since $0\in\DSS,$ we see that by the contraposition of $\ref{funzero},$ $f$ is 0 on $\co\pi(\FS).$ Since, by assumption, this hull in dense in $\DSS,$ any element in the range of $f$ over $\DSS$ is achieved as the limit of a net of $0,$ so must be $0$ on its entirety.
\end{proofm}

\begin{thm}
Let $\mathcal{T}$ be a locally convex topology on $B(H).$ If 
\[\oco^\mathcal{T}\pi(\FS)=\DSS,\]
then $\mathcal{T}$ is finer than ultraweak and coarser than strong$^*$ on $\DSS.$
\end{thm}

\begin{proofm}
Firstly, since $\mathcal{T}$ makes $B(H)$ Hausdorff, it must make the entry-wise functionals continuous, i.e. it must be finer than Kendall's topology, which we know to coincide with the weak/ultraweak on $\DSS.$
Now let $F=(B(H),\mathcal{T})^*.$ Notice, we know that $\mathcal{T}$ is coarser than $\tau(B(H),F),$ by the Mackey-Arens theorem. So take $\mathcal{T}=\tau(B(H),F)$ without loss of generality. In this topology, convergence is defined as follows: $u^{(\lambda)}\rightarrow u$ iff for every convex, $\sigma(F,B(H))-$compact $K\subseteq F,$ $\sup_{f\in K}|f(u^{(\lambda)}-u)|\rightarrow 0.$ However, by $\ref{closzero},$ for all $f\in K\setminus B(H)_*,$ $f(u^{(\lambda)}-u)=0$ for all $u,u^{(\lambda)}\in\DSS.$ Moreover, since $B(H)_*$ is itself convex and $\sigma(B(H)^*,B(H))-$closed, this convergence is equivalent to that of $\tau(B(H),F\cap B(H)_*)$ for nets in $\DSS.$ Therefore, $\mathcal{T}$ is coarser than the Mackey topology on $\DSS,$ which again coincides with the strong$^*$ topology.
\end{proofm}

Now we will describe the exposed points of $\DSS$ and $\DS$ in each of the topologies. It turns out that they coincide with the extremes of either in the weak topology, so the same holds for the rest of the topologies.

\begin{thm}
$\mathrm{exp}_w(\DSS)=\PS$ and $\mathrm{exp}_w(\DS)=\pi(\S)$
\end{thm}

\begin{proofm}
Let $u\in\PS.$ 
For any $J\subseteq\mathbb{N},$ define $x^J\in H_1$ by
\[x^J_j=\sum_{j\in J}2^{-\frac{j}{2}}e_j.\]
Let $I$ be the index underlying the support projection of $u.$ Now, define the weakly continuous functional by
\[f(v)=\langle vx^I,ux^I\rangle-\langle vx^{I^c},x^\mathbb{N}\rangle.\]
Our hypothesis is that $f$ reaches its unique maximum in $\DSS$ at $u.$ It suffices to show that this holds in $\PS,$ by theorem $\ref{coh}.$ Let $v\in \PS\setminus\{u\}.$
\\ Case 1: If $vx^I\neq ux^I,$ then by positive definiteness, 

$\langle vx^I,ux^I\rangle <\langle ux^I,ux^I\rangle.$ Of course, $\langle vx^{I^c},x^\mathbb{N}\rangle\geq 0,$ so $f(v)<f(u).$
\\ Case 2: If $vx^I=ux^I,$ then $vx^{I^c}\neq ux^{I^c}=0.$ Thus, for some $j\in\mathbb{N},$ $(vx^{I^c})_j>0$ and positive definiteness yields 
$(ux^\mathbb{N})_j>0,$ so $\langle vx^{I^c},x^\mathbb{N}\rangle>0.$ Again, $f(v)<f(u).$ Notice that the same must hold for $u\in\pi(\S)$ over $\DS.$
\end{proofm}

\section{The Affine Hull of $\S$}

\begin{defn}\label{matrx}
For $\Lambda\subseteq\mathbb{C},$ let $B_\Lambda(H)$ be the bounded operators on $H$ with all coefficients w.r.t $E$ in $\Lambda.$ Notice, in the case $\Lambda=\mathbb{R},$ this is a real C$^*$-algebra, as it can be $*$-isomorphically identified with $B(\tilde{H}),$ where $\tilde{H}=\overline{\mathrm{span}_\mathbb{R}}(E).$
\end{defn}

We can express $\PS$ in the terms set in definition $\ref{matrx},$ which elegantly characterizes them in the unit ball of $B(H).$

\begin{prop}
$\PS=B_{\{0,1\}}(H)_1.$
\end{prop}

\begin{proofm}
Clearly, $\PS\subseteq B_{\{0,1\}}(H)_1.$
For the other direction, let 
\\ $u\in B_{\{0,1\}}(H)_1.$ 
Then the sum of the squares of the rows and columns cannot exceed $1,$ so every row and column has at most one $1$ entry.
\end{proofm}

\begin{lem}\label{strdense}
Let $A$ be a real $*$-subalgebra of $B_\mathbb{R}(H).$ Then $A$ is strongly dense in its real double commutant $A''_\mathbb{R}=A''\cap B_\mathbb{R}(H).$
\end{lem}

\begin{proofm}
Since we know $B_\mathbb{R}(H)$ is real $*$-isomorphic to $B(\tilde{H}),$ it suffices to show that the strong operator topology induced by $H$ on $B_\mathbb{R}(H)$ coincides with strong operator topology induced by $\tilde{H}$ on $B(\tilde{H}).$ Better yet, it suffices to show that the latter is finer than the former. Let $a_\lambda,a\in B_\mathbb{R}(H)$ and suppose for all $y\in\tilde{H},$ $\lim_\lambda\|(a_\lambda-a)y\|=0.$
Let $x\in H.$
\begin{align*}
\|(a_\lambda-a)x\|^2 &= \sum_{k=1}^\infty|\langle (a_\lambda-a)x,e_k\rangle|^2 \\
&= \sum_{k=1}^\infty\left(|\langle (a_\lambda-a)\Re(x),e_k\rangle+i\langle (a_\lambda-a)\Im(x),e_k\rangle|^2\right)\\
&= \|(a_\lambda-a)\Re(x)\|^2+\|(a_\lambda-a)\Im(x)\|^2.
\end{align*}
But we know $\Re(x),\Im(x)\in \tilde{H},$ so $\lim_\lambda\|(a_\lambda-a)x\|=0.$ Now the classical result from C$^*$ theory follows. \cite[115-116]{MU90}
\end{proofm}

\begin{cor}\label{macdense}
Let $A$ be a real $*$-subalgebra of $B_\mathbb{R}(H).$ Then $A$ is Mackey dense in its real double commutant $A''_\mathbb{R}=A''\cap B_\mathbb{R}(H).$
\end{cor}

\begin{proofm}
Let $A_1=\{a\in A\mid \|a\|\leq 1\}.$ By the same alterations made in lemma $\ref{strdense},$ Kaplanksy's density theorem extends to real $*$-subalgebras of $B_\mathbb{R}(H).$ \cite[131]{MU90}
Then $A_1$ is strongly dense in $(A'')_1,$ and since $A_1$ is convex and bounded, its closure coincides with the Mackey closure, so $A_1$ is Mackey dense in $(A''_\mathbb{R})_1.$ \cite[70,153]{TA79}
\end{proofm}

\begin{prop}
$\oaff^{M}\pi(\FS)=B_\mathbb{R}(H)$
\end{prop}

\begin{proofm}
Notice that
\begin{align*}
\oaff^{M}\pi(\FS) \supseteq \oco^{M}\pi(\FS) &\supseteq B_{\{0,1\}}(H)_1,
\end{align*}
by prop. $\ref{hullc}.$ Thus, $0\in \oaff^{M}\pi(\FS),$ so the affine hull is the real span,
\[\oaff^{M}\pi(\FS)=\overline{\mathrm{span}}^{M}_\mathbb{R}(\pi(\FS))=\overline{\mathbb{R}[\pi(\S)]}^{M}.\]
Finally, $\mathbb{R}[\pi(\S)]$ is a real $*$-subalgebra of $B_\mathbb{R}(H),$ so by corollary $\ref{macdense},$ $\mathbb{R}[\pi(\S)]$ is Mackey dense in its real double commutant, 
\[\mathbb{R}[\pi(\S)]_\mathbb{R}''=\pi(\S)_\mathbb{R}''=\mathbb{R}\mathrm{Id}_H'=B_\mathbb{R}(H).\]
\end{proofm}

This fact is quite remarkable. In finite dimensions, the affine hull of the symmetric group is the smallest affine space containing the Birkhoff polytope $B_n$, $\aff\pi(S_n))=\aff(B_n),$ which famously has dimension $(n-1)^2$ embedded in $n^2$ real-dimensional space. Heuristically, the codimension of $\aff\pi(S_n)$ in its matrix space tends towards infinity in the limit of $n.$ However, for any topology which the dual space of $B(H)$ is contained in the predual $B(H)_*,$ the smallest closed affine space containing the direct limit of these polytopes has real codimension 0. Moreover, as we have shown, this space is a vector space only if the relative topology is coarser than that of the strong$^*$ on the doubly substochastics.

\section{Conclusion}
Initially motivated by a conditionally convergent trace on $B(H),$ this paper was an outgrowth of adjacent research. However, it quickly became apparent that the results herein were fruitful in their own right. We have demonstrated that Birkhoff's Problem 111 truly requires the use of a non-standard topology, such as that which Kendall gave, confirming Isbell's suspicion \cite{IS55}. Moreover, the use of the Arens-Mackey topology has shown how much room Kendall's secondary theorem had to improve. Not only did we find a finer topology for which his theorem holds, we found two, neither of which are compact. To this end, we have also shown that the larger of these is maximal. That is, any locally convex topology finer than Mackey must not be strictly finer on the doubly substochastics. In regards to the doubly substochastics themselves, we have seen that they affinely generate the entire (real) algebra, a fact which is nowhere near true of the topologies which only generate the doubly stochastics. Finally, this text has raised a few questions: 
\[\]
\\ 1. Does there exist a locally convex topology (coarser than the norm topology on and compatible with $B(H)$) on $\DSS$ which is both strictly finer than Mackey and has $\co\pi(\S)$ as a dense subset? 
\\ 2. Does there exist a locally convex topology such that $\oco\pi(\S)=\DS$?
\\ 3. Can $\ext(\DSS)\subseteq \PS$ be proven independent of the axiom of choice? 
\\ 4. Does the norm closure of $\co\pi(\S)$ coincide with Isbell's line-sum norm \cite{IS62}? If not, what is the additional condition on $\DS$ which defines it?

\end{document}